\ifx\documentclass\undefined
\documentstyle[12pt]{article}
\else
\documentclass[12pt]{article}
\usepackage{latexsym}
\usepackage{amsfonts}
\usepackage{amsmath}
\usepackage{amssymb}
\fi

\author{ K\'aroly Bezdek 
\thanks{Partially supported by a Natural Sciences and 
Engineering Research Council of Canada Discovery Grant}}

\font\tenBbb=msbm10 at 12pt         \font\sevenBbb=msbm9    \font\fiveBbb=msbm7
\newfam\Bbbfam
\textfont\Bbbfam=\tenBbb \scriptfont\Bbbfam=\sevenBbb
\scriptscriptfont\Bbbfam=\fiveBbb

\def\kkk{\null\hfill $\Box$\smallskip}

\newcommand{\proof}{{\noindent\bf Proof:{\ \ }}}

\newtheorem{theorem}{Theorem}[section]
\newtheorem{lemma}[theorem]{Lemma}

\newtheorem{cor}[theorem]{Corollary}
\newtheorem{con}[theorem]{Conjecture}
\newtheorem{prob}[theorem]{Problem}

\title{Illuminating spindle convex bodies and minimizing the volume of spherical sets of constant width
}

\begin{document}

\maketitle

\date

\begin{abstract}
A subset of the $d$-dimensional Euclidean space having nonempty interior is called a {\it spindle convex body} if it is the intersection of (finitely or infinitely many) congruent $d$-dimensional closed balls. The spindle convex body is called a ``fat''
one, if it contains the centers of its generating balls. The core part of this paper is an extension of Schramm's theorem and its proof on illuminating convex bodies of constant width to the family of ``fat'' spindle convex bodies.
\end{abstract}

\medskip

\section{Introduction}

Let $\mathbf{K}$ be a convex body (i.e., a compact convex set with nonempty interior) in the $d$-dimensional Euclidean space $\mathbb{E}^{d}$, $d\geq 2$. According to Boltyanski \cite{Bol1} the direction $\mathbf{v}\in \mathbb{S}^{d-1}$  (i.e., the unit vector $\mathbf{v}$ of $\mathbb{E}^{d}$) illuminates the boundary point $\mathbf{b}$ of $\mathbf{K}$ if the halfline emanating from $\mathbf{b}$ having direction vector $\mathbf{v}$ intersects the interior of $\mathbf{K}$, where $\mathbb{S}^{d-1}\subset\mathbb{E}^{d}$ denotes the $(d-1)$-dimensional unit sphere centered at the origin ${\bf o}$ of $\mathbb{E}^{d}$. Furthermore, the directions $\mathbf{v}_1, \mathbf{v}_2, \dots , \mathbf{v}_n$ illuminate $\mathbf{K}$ if each boundary point of $\mathbf{K}$ is illuminated by at least one of the directions $\mathbf{v}_1, \mathbf{v}_2, \dots , \mathbf{v}_n$. Finally, the smallest $n$ for which there exist $n$ directions that illuminate $\mathbf{K}$ is called the {\it illumination number} of $\mathbf{K}$ denoted by $I(\mathbf{K})$. An equivalent but somewhat different looking concept of illumination was introduced by Hadwiger in \cite{H60}. There he proposed to use point sources instead of directions for the illumination of convex bodies. Based on these circumstances the following conjecture, that was independently raised by Boltyanski \cite{Bol1} and Hadwiger \cite{H60} in 1960, is called the {\it Boltyanski--Hadwiger Illumination Conjecture}: The illumination number $I(\mathbf{K})$ of any convex body $\mathbf{K}$ in $\mathbb{E}^{d}$, is at most $2^d$ and $I(\mathbf{K})=2^d$ if and only if $\mathbf{K}$ is an affine $d$-cube.

Let $\mathbf{K}$ be a convex body in $\mathbb{E}^{d}$ and let $F$ be a face of
$\mathbf{K}$ (i.e., let $F$ be the intersection of a supporting hyperplane of $\mathbf{K}$ with the boundary of $\mathbf{K}$). Recall that the {\it Gauss image} $\nu (F)$ of the face $F$ is the set of all points (i.e. unit vectors) ${\bf u}\in \mathbb{S}^{d-1}\subset \mathbb{E}^{d}$ with the property that the supporting hyperplane of $\mathbf{K}$ with outer normal vector ${\bf u}$ contains $F$. It is easy to see that the Gauss images of distinct faces of $\mathbf{K}$ have disjoint relative interiors in $\mathbb{S}^{d-1}$ and $\nu (F)$ is compact and spherically convex for any face $F$. (Recall that a set $Y\subset \mathbb{S}^{d-1}$ is spherically convex if it is contained in an open hemisphere of $\mathbb{S}^{d-1}$ and for every $\mathbf{y}_1,\mathbf{y}_2\in Y$ the shorter great-circular arc of $\mathbb{S}^{d-1}$ connecting $\mathbf{y}_1$ with $\mathbf{y}_2$ is in $Y$.) Now, let $Y\subset \mathbb{S}^{d-1}$ be a set of finitely many points. Then the {\it covering radius} of $Y$ is the smallest positive real number $r$ with the property that the family of $(d-1)$-dimensional closed spherical balls of (angular) radii $r$ centered at the points of $Y$ cover $\mathbb{S}^{d-1}$. The following, rather basic principle, seems to be new and can be quite useful for estimating the illumination numbers of some convex bodies in particular, in low dimensions.

\begin{theorem}\label{36}
Let $\mathbf{K}\subset \mathbb{E}^{d}$, $d\geq 3$, be a convex body and let $r$ be a positive real number with the property that the Gauss image $\nu (F)$ of any face $F$ of $\mathbf{K}$ can be covered by a $(d-1)$-dimensional closed spherical ball of (angular) radius $r$ in $\mathbb{S}^{d-1}$. Moreover, assume that there exist $k$ points of $\mathbb{S}^{d-1}$ with covering radius $R$ satisfying the inequality $r+R\le \frac{\pi}{2}$. Then $I(\mathbf{K})\le k$.
\end{theorem}

In what follows we are going to study sets called {\it spindle convex bodies}. Based on the recent paper \cite{blnp} of the author, L\'angi, Nasz\'odi and Papez we can introduce them as follows. A subset of $\mathbb{E}^{d}$ having nonempty interior is called a {\it spindle convex body} if it is the intersection of (finitely or infinitely many) congruent $d$-dimensional closed balls. Here without loss of generality we assume that the congruent balls generating our spindle convex bodies are all of unit radii. Also, it is convenient to use the notation $\mathbf{B}[X]$ for the spindle convex body that is the intersection of the closed $d$-dimensional unit balls centered at the points of the compact set $X\subset \mathbb{E}^{d}$. For a comprehensive list of properties of spindle convex bodies we refer the interested reader to \cite{blnp}.

Now, let us take the spindle convex body $\mathbf{B}[X]$ in $\mathbb{E}^{3}$.  First, observe that if the Euclidean diameter ${\rm diam}(X)$ of $X$ satisfies the inequality ${\rm diam}(X)\le 0.577$ (resp., ${\rm diam}(X)\le 0.774$), then for the spherical diameter ${\rm Sdiam}(\nu (F))$ of the Gauss image $\nu (F)$ of an arbitrary
face $F$ of $\mathbf{B}[X]$ the inequality 
$${\rm Sdiam}(\nu (F))\le 2\arcsin(\frac{0.577}{2}) <33.5364^{\circ}$$ 
$${\rm (resp.,} \ {\rm Sdiam}(\nu (F))\le 2\arcsin(\frac{0.774}{2}) <45.5360^{\circ}{\rm  )}$$ holds. (We note that for the purpose of this discussion we use the degree measure for angles following \cite{FoTa}.) Thus, using the spherical Jung theorem \cite{Dekster}, we obtain that the Gauss image $\nu (F)$ of any face $F$ of $\mathbf{B}[X]$ can be covered by a $2$-dimensional closed spherical disk of (angular) radius $\le \arcsin\frac{0.577}{\sqrt{3}}<19.459^{\circ}$ (resp., $\le \arcsin\frac{0.774}{\sqrt{3}}<26.543^{\circ}$). Second, recall the well-known spherical codes (see \cite{FoTa}) according to which on $\mathbb{S}^{2}$ there are $4$ (resp., $5$) points with covering radius $<70.529^{\circ}$ (resp., $<63.435^{\circ}$). Hence, Theorem~\ref{36} leads us to the following statement. 

\begin{cor}\label{19-1}
Let $\mathbf{B}[X]$ be a spindle convex body in $\mathbb{E}^{3}$.   
\item(i) If $0<{\rm diam}(X)\le 0.577$, then $I(\mathbf{B}[X])=4$;
\item(ii) If $0.577<{\rm diam}(X)\le 0.774$, then $I(\mathbf{B}[X])\le 5$.
\end{cor}

The related statement that if $0<{\rm diam}(X)\le 1$, then $I(\mathbf{B}[X])\le 6$ has already been proved in \cite{blnp}. Clearly, Corollary~\ref{19-1} suggests to attack the Illumination Conjecture for spindle convex bodies in $\mathbb{E}^{3}$ by letting $0<{\rm diam}(X)<2$ to get arbitrarily close to $2$ while satisfying $0<{\rm cr}(X)<1$, where ${\rm cr}(X)$ denotes the radius of the unique smallest $3$-dimensional closed ball containing $X$. In connection with this, it is natural to expect that the illumination number of any spindle convex body in $\mathbb{E}^{3}$ is always strictly less than $8$. Moreover, for the sake of completeness we mention that the best known upper bound on the illumination numbers of $3$-dimensional convex bodies is due to Papadoperakis \cite{PA} stating that the illumination number of any convex body in $\mathbb{E}^{3}$ is at most $16$. This happens to be the best known upper bound for the illumination numbers of $3$-dimensional spindle convex bodies as well. For more information on the status of the Illumination Conjecture in $\mathbb{E}^{3}$ we refer the interested reader to \cite{Be06} and the relevant references listed there.

It is rather natural to expect that estimates similar to Corollary~\ref{19-1} exist in higher dimensions. For more details on that we refer the interested reader to the recent paper \cite{BeKi09} of the author and Kiss. However, the following approach is more efficient if the dimension is sufficiently large. Before stating our result, we briefly outline the status of the Illumination Conjecture in higher dimensions. (For a more complete picture on that we refer the interested reader to \cite{Be06} and the relevant references listed there.) The current best upper bound for the illumination numbers of convex bodies in higher dimensions has been obtained by Rogers using the main result of \cite{ER} combined with some observations from \cite{ER62} and with the inequality of Rogers and Shephard \cite{RS} on the volume of difference bodies, and reads as follows. If $\mathbf{K}$ is an arbitrary convex body in the $d$-dimensional Euclidean space $\mathbb{E}^{d}$, $d\geq 2$, then
$$I(\mathbf{K})\leq \frac{{\rm vol}_d(\mathbf{K}-\mathbf{K})}{{\rm vol}_d(\mathbf{K})}(d\ln d+d\ln\ln d+5d)\le\binom{2d}{d}(d\ln d+d\ln\ln d+5d).$$
Moreover, for sufficiently large $d$, $5d$ can be replaced by $4d$. We mention also the inequality $I(\mathbf{K})\leq(d+1)d^{d-1}-(d-1)(d-2)^{d-1}$ due to Lassak \cite{L88}, which is valid for an arbitrary convex body $\mathbf{K}$ in $\mathbb{E}^{d}$, $d\geq 2$. (Actually, Lassak's estimate is (somewhat) better than the estimate of Rogers for some small values of $d$). Note that, from the point of view of the Illumination Conjecture, the estimate of Rogers is nearly best possible for centrally symmetric convex bodies, since in that case $\frac{{\rm vol}_d(\mathbf{K}-\mathbf{K})}{{\rm vol}_d(\mathbf{K})}=2^d$. However, most convex bodies are far from being symmetric and so, in general, one may wonder whether the Illumination Conjecture is true at all, in particular, in high dimensions. Thus, it was an important progress, when Schramm \cite{Schr} managed to prove the Illumination Conjecture for all convex bodies of constant width in dimension greater than or equal to $16$. In fact, he has proved the following inequality. If $\mathbf{W}$ is an arbitrary convex body of constant width in $\mathbb{E}^{d}, d\ge 3$, then $$I(\mathbf{W})< 5d\sqrt{d}(4+\ln d)\left(\frac{3}{2}\right)^{\frac{d}{2}}.$$ 

By taking a closer look of the proof of the above upper bound of Schramm published in \cite{Schr}, and making the necessary modifications it turnes out that the estimate in question can be somewhat improved, but more importantly it can be extended to the family of ``fat'' spindle convex bodies, which is much larger than the family of convex bodies of constant width. Thus, we have the following theorem.

\begin{theorem}\label{19-2}
Let $\mathbf{B}[X]$ be an arbitrary spindle convex body in $\mathbb{E}^{d}$, $d\ge 3$, with ${\rm diam}(X)$ $\le 1$. Then $$I(\mathbf{B}[X])< 4\left(\frac{\pi}{3}\right)^{\frac{1}{2}}d^{\frac{3}{2}}(3+\ln d)\left(\frac{3}{2}\right)^{\frac{d}{2}}<5d^{\frac{3}{2}}(4+\ln d)\left(\frac{3}{2}\right)^{\frac{d}{2}}.$$
\end{theorem}

On the one hand, $4\left(\frac{\pi}{3}\right)^{\frac{1}{2}}d^{\frac{3}{2}}(3+\ln d)\left(\frac{3}{2}\right)^{\frac{d}{2}}<2^d$ for all $d\ge 15$. (Moreover, for every $\epsilon>0$ if $d$ is sufficiently large, then $I(\mathbf{B}[X])<\left(\sqrt{1.5}+\epsilon\right)^d=(1.224\ldots +\epsilon)^d$.) On the other hand, based on the elegant construction of Kahn and Kalai \cite{KaKa93}, it is known (see \cite{AiZi}), that if $d$ is sufficiently large, then there exists a finite subset $X''$ of $\{0,1\}^d$ in $\mathbb{E}^{d}$ such that any partition of $X''$ into parts of smaller diameter requires more than $(1.2)^{\sqrt{d}}$ parts. Let $X'$ be the (positive) homothetic copy of $X''$ having unit diameter and let $X$ be the (not necessarily unique) convex body of constant width one containing $X'$. Then it follows via standard arguments that $I(\mathbf{B}[X])>(1.2)^{\sqrt{d}}$ with $X=\mathbf{B}[X]$.

One of the key steps in the proof of Theorem~\ref{19-2}, presented in the relevant section of this paper, is Lemma~\ref{lower-bound-for-the-spherical-Kakeya-Pal-problem}. In fact, a better lower bound for Lemma~\ref{lower-bound-for-the-spherical-Kakeya-Pal-problem} could lead to an improvement in the exponential factor $\left(\frac{3}{2}\right)^{\frac{d}{2}}$ of Theorem~\ref{19-2}. As the underlying spherical geometry problem of Lemma~\ref{lower-bound-for-the-spherical-Kakeya-Pal-problem} might be of independent interest we phrase it in a slightly different but equivalent way and make some comments. In order to do so we recall some standard terminology. By a {\it convex body} $\mathbf{C}$ in $\mathbb{S}^{d-1}$  we understand the intersection $\mathbb{S}^{d-1}\cap\mathbf{C}_{\bf o}$, where $\mathbf{C}_{\bf o}$ stands for a line-free $d$-dimensional closed convex cone with apex ${\bf o}$ in $\mathbb{E}^{d}$. We denote by ${\cal K}_{\mathbb{S}^{d-1}}$ the space of all convex bodies in $\mathbb{S}^{d-1}$, equipped with the Hausdorff metric. $\mathbf{L}\subset \mathbb{S}^{d-1}$ is called a {\it lune} of $\mathbb{S}^{d-1}$ if it is the intersection of two (distinct) closed hemispheres of $\mathbb{S}^{d-1}$ having nonempty interior. The {\it width} of $\mathbf{L}$ is simply the angular measure of the dihedral angle ${\rm pos}(\mathbf{L})$, where ${\rm pos}(\cdot)$ refers to the positive hull of the corresponding set in $\mathbb{E}^{d}$. The {\it minimal width} ${\rm Swidth}(\mathbf{C})$ of $\mathbf{C}\in {\cal K}_{\mathbb{S}^{d-1}}$ is the smallest width of the lunes that contain $\mathbf{C}$. Also, we say that $\mathbf{C}\in {\cal K}_{\mathbb{S}^{d-1}}$ is a {\it convex body of constant width $w$} if $w={\rm Swidth}(\mathbf{C})={\rm Sdiam}(\mathbf{C})$, where ${\rm Sdiam}(\cdot)$ refers to the spherical diameter of the corresponding set in $\mathbb{S}^{d-1}$. For $\mathbf{C}\in {\cal K}_{\mathbb{S}^{d-1}}$ the {\it polar body} $\mathbf{C}^*$ of $\mathbf{C}$ is defined by 
$$
\mathbf{C}^*:=\{\mathbf{x}\in \mathbb{S}^{d-1}\ |\ \langle \mathbf{x}, \mathbf{c}\rangle\le 0\ {\rm for}\  {\rm all}\ \mathbf{c}\in \mathbf{C}\},
$$ 
where $\langle\cdot,\cdot\rangle$ refers to the canonical inner product in $\mathbb{E}^{d}$. (The induced canonical Euclidean norm on $\mathbb{E}^d$ will be denoted by $\|\cdot\|$.)
Clearly, $\mathbf{C}^*\in {\cal K}_{\mathbb{S}^{d-1}}$. Now, the problem studied in Lemma~\ref{lower-bound-for-the-spherical-Kakeya-Pal-problem} is equivalent to the following. (Actually, for a proof of the equivalence one can use the theorem proved in \cite{Dekster*} according to which any subset of $\mathbb{S}^{d-1}$ having spherical diameter $0<w\le\frac{\pi}{2}$ can be covered by a convex body of constant width $w$ in $\mathbb{S}^{d-1}$. Moreover, the polar body of such a convex body is of constant width $\frac{\pi}{2}\le\pi-w<\pi$.) Let the positive real $\frac{\pi}{2}\le w^*<\pi$ and the positive integer $d\ge 3$ be given. Then find the minimum volume convex body of constant width $w^*$ in $\mathbb{S}^{d-1}$. In fact, the question makes sense to ask for all $0<w^*<\pi$. Thus, we have arrived at the following quite basic volume problem, whose Euclidean counterpart has been much better studied and is also better known (see for example \cite{BLO07}).

\begin{prob}\label{Blaschke's-spherical-problem}
For $0< w^*<\pi$ and $d\ge 3$ find the minimum volume convex body of constant width $w^*$ in $\mathbb{S}^{d-1}$.
\end{prob}

Problem~\ref{Blaschke's-spherical-problem} has been solved by Blaschke on $\mathbb{S}^{2}$ (i.e., for $d=3$) in \cite{Bl1915}. Blaschke's theorem (\cite{Bl1915}) can be summarized as follows. Among all convex domains of constant width $0<w^*\le\frac{\pi}{2}$ on $\mathbb{S}^{2}$, the Reuleaux triangle of constant width $w^*$ has the smallest area. Moreover, among all convex domains of constant width $\frac{\pi}{2}<w^*<\pi$ of $\mathbb{S}^{2}$ the smallest area belongs to the one which is obtained as the outer parallel domain of radius $w^*-\frac{\pi}{2}$ of the Reuleaux triangle of width $\pi-w^*$. We have the following partial extension of this theorem of Blaschke to $\mathbb{S}^{3}$.

\begin{theorem}\label{19-3}
Let $0<w<\pi$ be given. Then the volume of the convex body $\mathbf{C}$ of constant width $w$ in $\mathbb{S}^{3}$ is minimal among all convex bodies of constant width $w$ of $\mathbb{S}^{3}$ if and only if the polar body $\mathbf{C}^*$ of constant width $\pi-w$ has minimal volume among all convex bodies of constant width $\pi-w$ of $\mathbb{S}^{3}$. 
\end{theorem}

Thus, Theorem~\ref{19-3} implies that if $d=4$, then it is sufficient to investigate Problem~\ref{Blaschke's-spherical-problem} for convex bodies of constant widths $0<w\le\frac{\pi}{2}$ in $\mathbb{S}^{3}$. The question whether a statement similar to Theorem~\ref{19-3} holds in spherical spaces of dimensions $4$ and higher remains open. Finally, we wish to call the reader's attention to the following special case of Problem~\ref{Blaschke's-spherical-problem} that is strikingly simple to phrase in all spherical dimensions.

\begin{con}\label{Bezdek's-Spherical-Volume-Conjecture}
Among all convex bodies of constant width $\frac{\pi}{2}$ in $\mathbb{S}^{d-1}$, $d\ge 4$, the $(d-1)$-dimensional regular simplex of edge length $\frac{\pi}{2}$ has the smallest volume. 
\end{con}

\section{Proof of Theorem 1.1}

Recall the following well-known observation on illumination. For the convenience of the reader and for notational reasons we include its short proof. (For more information on different approaches to illumination we refer the interested reader to \cite{Be06} and the relevant references listed there.) Let the open ball centered at the point $\mathbf{p}\in\mathbb{S}^{d-1}$ having (angular) radius $0<\alpha<\pi$ in the spherical space $\mathbb{S}^{d-1}$ be denoted by $C(\mathbf{p}, \alpha)$ and let us call it the {\it open spherical cap} of $\mathbb{S}^{d-1}$ with center $\mathbf{p}$ and radius $\alpha$. In particular, $C(\mathbf{p}, \frac{\pi}{2})$ will be called the
{\it open hemisphere} of $\mathbb{S}^{d-1}$ with center $\mathbf{p}$.

\begin{lemma}\label{separation-for-Gauss-images}
Let $\mathbf{K}$ be a convex body in $\mathbb{E}^{d}$, $d\ge 3$, and let $\mathbf{b}\in {\rm bd}(\mathbf{K})$ be an arbitrary boundary point of $\mathbf{K}$. Moreover, let $F_{\mathbf{b}}$ denote the smallest dimensional face of $\mathbf{K}$ containing $\mathbf{b}$. Then $\mathbf{b}$ is illuminated by the direction $\mathbf{v}\in \mathbb{S}^{d-1}$ if and only if 
$$
\nu \left(F_{\mathbf{b}}\right)\subset C\left(-\mathbf{v}, \frac{\pi}{2}\right).
$$
Furthermore, $I(\mathbf{K})$ is the smallest number of open hemispheres of $\mathbb{S}^{d-1}$ with the property that the Gauss image of each face of $\mathbf{K}$ is contained in at least one of the given open hemispheres.
\end{lemma}

\proof
It is convenient to use the following notation. For a set $A\subset \mathbb{S}^{d-1}$ let $A^+:=\{\mathbf{x}\in \mathbb{S}^{d-1}\ |\ \langle \mathbf{x}, \mathbf{y}\rangle>0 \ {\rm for \ all}\ \mathbf{y}\in A\}$.  

First, observe that the halfline emanating from $\mathbf{b}\in{\rm relint} (F_{\mathbf{b}})$ (with ${\rm relint} (\cdot)$ standing for the relative interior of the corresponding set) having direction vector $\mathbf{v}$ intersects the interior of $\mathbf{K}$ if and only if $-\mathbf{v}\in \nu \left(F_{\mathbf{b}}\right)^+$. Second, observe that $-\mathbf{v}\in \nu \left(F_{\mathbf{b}}\right)^+$ if and only if $\nu \left(F_{\mathbf{b}}\right)\subset C\left(-\mathbf{v}, \frac{\pi}{2}\right)$. This completes the proof of Lemma~\ref{separation-for-Gauss-images}.
\kkk

Now, we turn to the proof of Theorem~\ref{36}. Let $\{\mathbf{p}_1,\dots, \mathbf{p}_k\}$ be the family of points in $\mathbb{S}^{d-1}$ with covering radius $R$. Moreover, let $B_i\subset \mathbb{S}^{d-1}$  be the $(d-1)$-dimensional closed spherical ball of radius $R$ centered at the point $\mathbf{p}_i$ in $\mathbb{S}^{d-1}$, $1\le i\le k$. Finally, let $C_i$ be the open hemisphere of $\mathbb{S}^{d-1}$ with center $\mathbf{p}_i$, $1\le i\le k$. Based on Lemma~\ref{separation-for-Gauss-images} it is sufficient to show that the
Gauss image of each face of $\mathbf{K}$ is contained in at least one of the open hemispheres $C_i, 1\le i\le k$. 

Now, let $F$ be an arbitrary face of the convex body $\mathbf{K}\subset \mathbb{E}^{d}$, $d\geq 3$, and let $B_F$ denote the smallest $(d-1)$-dimensional closed spherical ball of $\mathbb{S}^{d-1}$ with center $\mathbf{f}\in \mathbb{S}^{d-1}$ which contains the Gauss image $\nu (F)$ of $F$. By assumption the radius of $B_F$ is at most $r$. As the family $\{B_i, 1\le i\le k\}$ of balls forms a covering of $\mathbb{S}^{d-1}$ therefore $\mathbf{f}\in B_j$ for some $1\le j\le k$. If in addition, we have that $\mathbf{f}\in {\rm Sint} (B_j)$ (where ${\rm Sint}(\cdot )$ denotes the (spherical) interior of the corresponding set in $\mathbb{S}^{d-1}$), then the inequality $r+R\le \frac{\pi}{2}$ implies that $\nu (F)\subset C_j$. If $\mathbf{f}$ does not belong to the interior of any of the sets $B_i, 1\le i\le k$, then clearly $\mathbf{f}$ must be a boundary point of at least $d$ sets of the family $\{B_i, 1\le i\le k\}$. Then either we find a $C_i$ containing $\nu (F)$ or we end up with $d$ members of the family $\{C_i, 1\le i\le k\}$ each being tangent to $B_F$ at some point of $\nu (F)$. Clearly, the later case can occur only for finitely many $\nu (F)$'s and so, by taking a proper congruent copy of the open hemispheres $\{C_i, 1\le i\le k\}$ within $\mathbb{S}^{d-1}$ (under which we mean to avoid finitely many so-called prohibited positions) we get that each 
$\nu (F)$ is contained in at least one member of the family $\{C_i, 1\le i\le k\}$. This completes the proof of Theorem~\ref{36}.
\kkk

\section{Proof of Theorem 1.3}

\subsection{On the boundary of spindle convex hulls}

Let $X\subset \mathbb{E}^d$, $d\ge 3$, be a compact set with ${\rm cr}(X)< 1$, where ${\rm cr}(X)$ denotes the radius of the smallest $d$-dimensional closed Euclidean ball containing $X$. For the following investigations it will be more proper to use the normal images than the Gauss images of the boundary points of $\mathbf{B}[X]$. The {\it normal image} $N_{\mathbf{B}[X]}(\mathbf{b})$ of an arbitrary boundary point $\mathbf{b}\in {\rm bd}\left(\mathbf{B}[X]\right)$ of $\mathbf{B}[X]$ is defined as
$$
N_{\mathbf{B}[X]}(\mathbf{b}):=-\nu(\{\mathbf{b}\})
$$
In other words, $N_{\mathbf{B}[X]}(\mathbf{b})\subset \mathbb{S}^{d-1}$ is the set of inward unit normal vectors of all hyperplanes that support $\mathbf{B}[X]$ at $\mathbf{b}$. Clearly, $N_{\mathbf{B}[X]}(\mathbf{b})$ is a closed spherically convex subset of $\mathbb{S}^{d-1}$. Moreover, Lemma~\ref{separation-for-Gauss-images} implies in a straighforward way that the direction $\mathbf{u}\in \mathbb{S}^{d-1}$ illuminates the boundary point $\mathbf{b}$ of the convex body $\mathbf{B}[X]$ if and only if $\mathbf{u}\in N_{\mathbf{B}[X]}(\mathbf{b})^+$.

We will need the following definitions and lemma from \cite{blnp}. Let ${\bf a}$ and ${\bf b}$ be two points in $\mathbb{E}^{d}$. If $\|{\bf a}-{\bf b}\|<2$, then the {\it (closed) spindle} of ${\bf a}$ and ${\bf b}$, denoted by $[{\bf a},{\bf b}]_s$, is defined as the union of circular arcs  with endpoints ${\bf a}$ and ${\bf b}$ which have radii at least one and are shorter than a semicircle. If $\|{\bf a}-{\bf b}\|=2$, then $[{\bf a},{\bf b}]_s:=\mathbf{B}^d[{\bf \frac{a+b}{2}}, 1]$, where $\mathbf{B}^d[{\bf p}, r]$ denotes the (closed) $d$-dimensional ball centered at ${\bf p}$ with radius $r$ in $\mathbb{E}^{d}$. If $\|{\bf a}-{\bf b}\|>2$, then we define $[{\bf a},{\bf b}]_s$ to be $\mathbb{E}^{d}$. Next, a set $\mathbf{C}\subset\mathbb{E}^{d}$ is called {\it spindle convex} if, for any pair of points ${\bf a}, {\bf b}\in \mathbf{C}$, we have that $[{\bf a}, {\bf b}]_s\subset \mathbf{C}$. Finally, let $X$ be a set in $\mathbb{E}^{d}$. Then the {\it spindle convex hull} of $X$ is the set defined by ${\rm conv}_s X := \bigcap \{ C\subset \mathbb{E}^{d} | X \subset C \ {\rm and}\ C \ {\rm is \ spindle \  convex \ in}\  \mathbb{E}^{d}\}$. Also, recall that $S^{d-1}(\mathbf{c},r)\subset\mathbb{E}^d$ denotes the $(d-1)$-dimensional sphere centered at $\mathbf{c}$ having radius $r$. A set $Y\subset S^{d-1}(\mathbf{c},r)$ is {\it spherically convex} if it is contained in an open hemisphere of $S^{d-1}(\mathbf{c},r)$ and for every $\mathbf{y}_1,\mathbf{y}_2\in Y$ the shorter great-circular arc of $S^{d-1}(\mathbf{c},r)$ connecting $\mathbf{y}_1$ with $\mathbf{y}_2$ is in $Y$. The {\it spherical convex hull} of a set $Y\subset S^{d-1}(\mathbf{c},r)$ is defined in the natural way and it
exists if, and only if, $Y$ is in an open hemisphere of $S^{d-1}(\mathbf{c},r)$. We denote it by ${\rm Sconv}(Y, S^{d-1}(\mathbf{c},r))$. The following lemma proved in \cite{blnp} describes some properties of the boundary of spindle convex hulls.

\begin{lemma}\label{lem:sconv}
Let $X\subset \mathbb{E}^d$ be a compact set. If ${\rm cr}(X)<1$ and $\mathbf{B}^d[\mathbf{q}, 1]$ is a closed unit ball containing $X$, then

\item(i) $X\cap S^{d-1}(\mathbf{q},1)$ is contained in an open hemisphere of $S^{d-1}(\mathbf{q},1)$ and

\item(ii) ${\rm conv}_s(X)\cap S^{d-1}(\mathbf{q},1)={\rm Sconv}(X\cap S^{d-1}(\mathbf{q},1), S^{d-1}(\mathbf{q},1)).$

\end{lemma}

Now, we are ready to prove the main lemma of this section.

\begin{lemma}\label{basic-illumination-spindle-hull}
Let $X\subset \mathbb{E}^d$, $d\ge 3$, be a compact set with ${\rm cr}(X)< 1$. Then the boundary of the spindle convex hull of $X$ can be generated as follows: $${\rm bd}\left({\rm conv}_s(X)\right)=\bigcup_{\mathbf{b}\in {\rm bd}\left(\mathbf{B}[X]\right)}\{ \mathbf{b}+\mathbf{y}\ |\ \mathbf{y}\in N_{\mathbf{B}[X]}(\mathbf{b})\}.$$
\end{lemma}

\proof
Let $\mathbf{b}\in {\rm bd}\left(\mathbf{B}[X]\right)$. Then (ii) of Lemma~\ref{lem:sconv} implies that
$$ \mathbf{b}+N_{\mathbf{B}[X]}(\mathbf{b})={\rm Sconv}(X\cap S^{d-1}(\mathbf{b},1), S^{d-1}(\mathbf{b},1))={\rm conv}_s(X)\cap S^{d-1}(\mathbf{b},1) .$$

This together with the fact that 
$$
\bigcup_{\mathbf{b}\in {\rm bd}\left(\mathbf{B}[X]\right)}N_{\mathbf{B}[X]}(\mathbf{b})=\mathbb{S}^{d-1}
$$
finishes the proof of Lemma~\ref{basic-illumination-spindle-hull}. 
\kkk

\subsection{On the Euclidean diameter of spindle convex hulls and normal images}

\begin{lemma}\label{diameter-of-spindle-hull}
If $X\subset \mathbb{E}^d$, $d\ge 3$, is a compact set with ${\rm diam}(X)\le 1$, then
$${\rm diam}\left({\rm conv}_s(X)\right)\le 1.$$
\end{lemma}
\proof
By assumption ${\rm diam}(X)\le 1$. Recall that Meissner \cite{Me11} has called a compact set $M\subset\mathbb{E}^d$ {\it complete} if ${\rm diam}(M\cup\{\mathbf{p}\})> {\rm diam}(M)$ for any $\mathbf{p}\in\mathbb{E}^d\setminus M$. He has proved in \cite{Me11} that any set of diameter $1$ is contained in a complete set of diameter $1$. Moreover, he has shown in \cite{Me11} that a compact set of diameter $1$ in $\mathbb{E}^d$ is complete if and only if it is of constant width $1$. These facts together with the easy observation that any convex body of constant width $1$ in $\mathbb{E}^d$ is in fact, a spindle convex set, imply that $X$ is contained in a convex body of constant width $1$ and any such convex body must necessarily contain ${\rm conv}_s(X)$. Thus, indeed ${\rm diam}\left({\rm conv}_s(X)\right)\le 1$. 
\kkk

For an arbitrary nonempty subset $A$ of $\mathbb{S}^{d-1}$ let
$$U_{\mathbf{B}[X]}(A):=\left(\bigcup_{N_{\mathbf{B}[X]}(\mathbf{b})\cap A\neq\emptyset} N_{\mathbf{B}[X]}(\mathbf{b})\right) \subset \mathbb{S}^{d-1}.$$

\begin{lemma}\label{diameter-of-union-of-normal-images}
Let $X\subset \mathbb{E}^d$, $d\ge 3$, be a compact set with ${\rm diam}(X)\le 1$ and 
let $\emptyset\neq A\subset\mathbb{S}^{d-1}$ be given. Then  
$${\rm diam} \left(U_{\mathbf{B}[X]}(A)\right)\le 1+{\rm diam} (A).$$
\end{lemma}

\proof
Let $\mathbf{y}_1\in N_{\mathbf{B}[X]}(\mathbf{b}_1)$ and $\mathbf{y}_2\in N_{\mathbf{B}[X]}(\mathbf{b}_2)$ be two arbitrary points of $U_{\mathbf{B}[X]}(A)$ with
$\mathbf{b}_1, \mathbf{b}_2 \in {\rm bd}\left(\mathbf{B}[X]\right)$. We need to show that $\|\mathbf{y}_1-\mathbf{y}_2\|\le  1+{\rm diam} (A)$. 

By Lemma~\ref{basic-illumination-spindle-hull} and by Lemma~\ref{diameter-of-spindle-hull} we get that $$\|(\mathbf{y}_1-\mathbf{y}_2)+(\mathbf{b}_1-\mathbf{b}_2) \|=  \|(\mathbf{b}_1+\mathbf{y}_1)-(\mathbf{b}_2+\mathbf{y}_2)\|\le 1.$$
Thus, the reverse triangle inequality yields that
$$\|\mathbf{y}_1-\mathbf{y}_2\|\le 1+ \|\mathbf{b}_2-\mathbf{b}_1 \|.$$
This means that in order to finish the proof of Lemma~\ref{diameter-of-union-of-normal-images} it is sufficient to show that $\|\mathbf{b}_2-\mathbf{b}_1 \|\le {\rm diam} (A)$. This can be done as follows. First,
note that the sets $\mathbf{b}_1 + N_{\mathbf{B}[X]}(\mathbf{b}_1)\subset {\rm bd}\left({\rm conv}_s(X)\right)$ and $\mathbf{b}_2 + N_{\mathbf{B}[X]}(\mathbf{b}_2)\subset {\rm bd}\left({\rm conv}_s(X)\right)$ are separated by the hyperplane $H$ of $\mathbb{E}^d$ that bisects the line segment connecting $\mathbf{b}_1$ to $\mathbf{b}_2$ and is perpendicular to it with $\mathbf{b}_1 +N_{\mathbf{B}[X]}(\mathbf{b}_1)$ (resp., $\mathbf{b}_2 + N_{\mathbf{B}[X]}(\mathbf{b}_2)$) lying on the same side of $H$ as $\mathbf{b}_2$ (resp., $\mathbf{b}_1$). (All this follows in a direct way from the observation that a unit ball centered at an arbitrary point of ${\rm conv}_s(X)$ contains $\mathbf{B}[X]$.) Second, assume that $\|\mathbf{b}_2-\mathbf{b}_1 \| > {\rm diam} (A)$. Then this assumption together with the separating hyperplane $H$ clearly imply that the Euclidean distance between the sets $N_{\mathbf{B}[X]}(\mathbf{b}_1)$ and $N_{\mathbf{B}[X]}(\mathbf{b}_2)$ is at least $\|\mathbf{b}_2-\mathbf{b}_1 \| > {\rm diam} (A)$, a contradiction (since by the assumption of Lemma~\ref{diameter-of-union-of-normal-images} we have that $N_{\mathbf{B}[X]}(\mathbf{b}_1)\cap A\neq\emptyset$ and $N_{\mathbf{B}[X]}(\mathbf{b}_2)\cap A\neq\emptyset$). This completes the proof of Lemma~\ref{diameter-of-union-of-normal-images}.
\kkk

\subsection{An upper bound for the illumination number}

Let $\mu_{d-1}$ denote the standard probability measure on $\mathbb{S}^{d-1}$ and define
$$V_{d-1} (t):=\inf\{ \mu_{d-1}(A^+)\ |\ A\subset \mathbb{S}^{d-1}, {\rm diam}(A)\le t\},$$
where $0<t\le\sqrt{2}$. Moreover, let $n_{d-1}(\epsilon)$ denote the minimum number of closed spherical caps of $\mathbb{S}^{d-1}$ having Euclidean diameter $\epsilon$ such that they cover $\mathbb{S}^{d-1}$, where $0<\epsilon\le 2$.

\begin{lemma}\label{upper-estimate-for-the-illumination-number}
$$ I(\mathbf{B}[X])\le 1+\frac{\ln \left(n_{d-1}(\epsilon)\right)}{-\ln \left(1-V_{d-1} (1+\epsilon)\right)}$$
holds for all $0<\epsilon \le \sqrt{2}-1$ and $d\ge 3$.
\end{lemma}

\proof
Let $\emptyset\ne A\subset \mathbb{S}^{d-1}$ be given with ${\rm diam}(A) \le 1+\epsilon \le \sqrt{2}$. Then the spherical Jung theorem \cite{Dekster} implies that $A$ is contained in a closed spherical cap of $\mathbb{S}^{d-1}$ having angular radius $0<\arcsin\sqrt{\frac{d-1}{d}}$ $<\frac{\pi}{2}$. Thus, $A^+$ contains a spherical cap of $\mathbb{S}^{d-1}$ having angular radius $\frac{\pi}{2}-\arcsin\sqrt{\frac{d-1}{d}}>0$ and of course, $A^+$ is contained in an open hemisphere of $\mathbb{S}^{d-1}$. Hence, $0<V_{d-1} (1+\epsilon)<\frac{1}{2}$ and so, the expression on the right in Lemma~\ref{upper-estimate-for-the-illumination-number} is well-defined.

Let $m$ be a positive integer satisfying
$$m> \frac{\ln \left(n_{d-1}(\epsilon)\right)}{-\ln \left(1-V_{d-1} (1+\epsilon)\right)} .$$
It is sufficient to show that $m$ directions can illuminate $\mathbf{B}[X]$. Let $n:=n_{d-1}(\epsilon)$ and let $A_1, A_2, \dots, A_n$ be closed spherical caps of $\mathbb{S}^{d-1}$ having Euclidean diameter $\epsilon$ and covering $\mathbb{S}^{d-1}$. By Lemma~\ref{diameter-of-union-of-normal-images} we have 
$${\rm diam} \left(U_{\mathbf{B}[X]}(A_i)\right)\le 1+\epsilon$$ 
for all $1\le i\le n$ and therefore
$$\mu_{d-1}\left(U_{\mathbf{B}[X]}(A_i)^+  \right)\ge V_{d-1} (1+\epsilon)$$
for all $1\le i\le n$. Let the directions $\mathbf{u}_1, \mathbf{u}_2, \dots , \mathbf{u}_m$ be chosen at random, uniformly and independently distributed on $\mathbb{S}^{d-1}$. Thus, the probability that $\mathbf{u}_j$ lies in $U_{\mathbf{B}[X]}(A_i)^+$ is equal to $\mu_{d-1}\left(U_{\mathbf{B}[X]}(A_i)^+ \right)\ge V_{d-1} (1+\epsilon)$. Therefore the probabilty that $U_{\mathbf{B}[X]}(A_i)^+$ contains none of the points $\mathbf{u}_1, \mathbf{u}_2, \dots , \mathbf{u}_m$ is at most $\left(1-V_{d-1} (1+\epsilon)\right)^m$. Hence, the probability $p$ that at least one $U_{\mathbf{B}[X]}(A_i)^+$ contains none of the points $\mathbf{u}_1, \mathbf{u}_2, \dots , \mathbf{u}_m$ satisfies
$$p\le \sum_{i=1}^n \left(1-V_{d-1} (1+\epsilon)\right)^m< n\left(1-V_{d-1} (1+\epsilon)\right)^{\frac{\ln \left(n\right)}{-\ln \left(1-V_{d-1} (1+\epsilon)\right)}}=1.$$
This shows that one can choose $m$ directions say, $\{\mathbf{v}_1, \mathbf{v}_2, \dots , \mathbf{v}_m\}\subset \mathbb{S}^{d-1}$, such that each set $U_{\mathbf{B}[X]}(A_i)^+$, $1\le i\le n$, contains at least one of them. We claim that the directions $\mathbf{v}_1, \mathbf{v}_2, \dots , \mathbf{v}_m$ illuminate $\mathbf{B}[X]$. Indeed, let $\mathbf{b}\in {\rm bd}\left( \mathbf{B}[X]\right)$. We will show that at least one of the directions $\mathbf{v}_1, \mathbf{v}_2, \dots , \mathbf{v}_m$ illuminates the boundary point $\mathbf{b}$. As the spherical caps $A_1, A_2, \dots, A_n$ form a covering of $\mathbb{S}^{d-1}$ therefore there exists an $A_i$ with $A_i\cap N_{\mathbf{B}[X]}(\mathbf{b})\neq\emptyset$. Thus, by definition $ N_{\mathbf{B}[X]}(\mathbf{b})\subset U_{\mathbf{B}[X]}(A_i)$ and therefore
$$ N_{\mathbf{B}[X]}(\mathbf{b})^+\supset U_{\mathbf{B}[X]}(A_i)^+ .$$
$U_{\mathbf{B}[X]}(A_i)^+$ contains at least one of the directions $\mathbf{v}_1, \mathbf{v}_2, \dots , \mathbf{v}_m$, say $\mathbf{v}_k$. Hence,
$$\mathbf{v}_k\in U_{\mathbf{B}[X]}(A_i)^+\subset N_{\mathbf{B}[X]}(\mathbf{b})^+$$
and so, Lemma~\ref{separation-for-Gauss-images} yields that indeed, $\mathbf{v}_k$ illuminates the boundary point $\mathbf{b}$ of $\mathbf{B}[X]$, finishing the proof of Lemma ~\ref{upper-estimate-for-the-illumination-number}.
\kkk

\subsection{Schramm's lower bound for the proper measure of polars of sets of given diameter in spherical space}

We need the following notation for the next statement. For $\mathbf{u}\in \mathbb{S}^{d-1}$ let $R_{\mathbf{u}}: \mathbb{E}^d\to \mathbb{E}^d$ denote the reflection about the line passing through the points $\mathbf{u}$ and $-\mathbf{u}$. Clearly, $R_{\mathbf{u}}(\mathbf{x})=2\langle\mathbf{x},\mathbf{u}\rangle\mathbf{u}- \mathbf{x}$ for all $\mathbf{x}\in \mathbb{E}^d$. As the following two lemmas are
taken from \cite{Schr} with some minor changes in notation we quote them without proof.

\begin{lemma}\label{reflection-principle}
Let $A\subset \mathbb{S}^{d-1}$ be a set of Euclidean diameter $0<{\rm diam}(A)\le t$ contained in the closed spherical cap $C[\mathbf{u}, \arccos a]\subset \mathbb{S}^{d-1}$ centered at $\mathbf{u}\in \mathbb{S}^{d-1}$ having angular radius $0<\arccos a<\frac{\pi}{2}$ with $0<a<1$ and $0<t\le 2\sqrt{1-a^2}$. Then
$$ A^+\cup R_{\mathbf{u}}(A^+)\supset C\left(\mathbf{u},\arctan\left(\frac{2a}{t}\right)\right)  .$$
\end{lemma}

\begin{lemma}\label{lower-bound-for-the-spherical-Kakeya-Pal-problem}
$$V_{d-1} (t)\ge\frac{1}{\sqrt{8\pi d}}\left(\frac{3}{2}+\frac{\left(2-\frac{1}{d} \right)t^2-2}{4-\left(2-\frac{2}{d} \right)t^2}\right)^{-\frac{d-1}{2}}$$
for all $0<t\le\sqrt{2}$ and $d\ge 3$.
\end{lemma}

\subsection{An upper bound for the number of sets of given diameter that are needed to cover spherical space}

The following (simple) estimate is well-known (see for example \cite{Schr}). We refer the interested reader for a proof to the proper section in \cite{Schr}.

\begin{lemma}\label{maximal-packing-estimate}
$$n_{d-1}(\epsilon)< \left(1+\frac{4}{\epsilon}\right)^d$$
for all $0<\epsilon \le 2$ and $d\ge 3$.
\end{lemma}

Actually, using \cite{Du07}, one can replace the inequality of Lemma~\ref{maximal-packing-estimate} by the stronger inequality $n_{d-1}(\epsilon)\le (\frac{1}{2} + {\rm o}(1))d\ln d\left(\frac{2}{\epsilon}\right)^d$. As this improves the estimate of Theorem ~\ref{19-2} only in a rather insignificant way, we do not introduce it here.

\subsection{The final upper bound for the illumination number}

Now, we are ready for the proof of Theorem ~\ref{19-2}. As $x<-\ln (1-x)$ holds for all $0<x<1$, therefore by Lemma~\ref{upper-estimate-for-the-illumination-number} we get that
$$ I(\mathbf{B}[X])\le 1+\frac{\ln \left(n_{d-1}(\epsilon)\right)}{-\ln \left(1-V_{d-1} (1+\epsilon)\right)}<1+\frac{\ln \left(n_{d-1}(\epsilon)\right)}{V_{d-1} (1+\epsilon)}$$
holds for all $0<\epsilon\le\sqrt{2}-1$ and $d\ge 3$. Now, let $\epsilon_0=\sqrt{\frac{2d}{2d-1}}-1$. As $0<\epsilon_0<\sqrt{2}-1$ holds for all $d\ge 3$, therefore Lemma~\ref{lower-bound-for-the-spherical-Kakeya-Pal-problem} and Lemma~\ref{maximal-packing-estimate} together with the easy inequality $\epsilon_0>\frac{4}{16d-1}$ yield that
$$I(\mathbf{B}[X])<1+\sqrt{8\pi d}\left(\frac{3}{2}\right)^{\frac{d-1}{2}}\ln \left(n_{d-1}(\epsilon_0)\right)$$ 
$$<1+\sqrt{8\pi d}\left(\frac{3}{2}\right)^{\frac{d-1}{2}}\ln \left( \left(1+\frac{4}{\epsilon_0}\right)^d\right)<1+\sqrt{8\pi d}\left(\frac{3}{2}\right)^{\frac{d-1}{2}}\ln \left( (16d)^d\right)$$
$$=1+4\sqrt{\frac{\pi}{3}}d\sqrt{d}\left(\frac{3}{2}\right)^{\frac{d}{2}}(\ln 16+\ln d)<
4\left(\frac{\pi}{3}\right)^{\frac{1}{2}}d^{\frac{3}{2}}(3+\ln d)\left(\frac{3}{2}\right)^{\frac{d}{2}},$$
finishing the proof of Theorem~\ref{19-2}.

\section{Proof of Theorem 1.5}

Let $\mathbf{C}\in {\cal K}_{\mathbb{S}^{3}}$ be an arbitrary convex body of constant width $0<w<\pi$ with sufficiently smooth boundary in $\mathbb{S}^{3}$. On the one hand, according to a classical result of Blaschke \cite{Bl1915} we have that
$$
M_1(\mathbf{C})+2V(\mathbf{C})=2\pi w,
$$
where $M_1(\mathbf{C})$ is the integral of the mean curvature (evaluated over the boundary of $\mathbf{C}$) and $V(\mathbf{C})$ denotes the (spherical) volume of $\mathbf{C}$. (See also formula (5.7) in \cite{Sa50}.) On the other hand, Allendoerfer \cite{Al48} has proved that (for not only the above $\mathbf{C}$, but actually, for any $\mathbf{C}\in {\cal K}_{\mathbb{S}^{3}}$ with sufficiently smooth boundary) we have also
$$
M_1(\mathbf{C})+V(\mathbf{C})+V(\mathbf{C}^*)=\pi^2.
$$
(See also formula (17.31) in \cite{Sa76}.) Clearly, the above two equations imply that
$$
\pi^2+V(\mathbf{C})=2\pi w+V(\mathbf{C}^*)
$$ 
holds for any convex body $\mathbf{C}$ of constant width $0<w<\pi$ in $\mathbb{S}^{3}$, from which Theorem~\ref{19-3} follows in a straightforward way.

\vspace{1cm}

\medskip

\noindent
K\'aroly Bezdek
\newline
Department of Mathematics and Statistics, University of Calgary, Canada,
\newline
Department of Mathematics, University of Pannonia, Veszpr\'em, Hungary,
\newline
and
\newline
Institute of Mathematics, E\"otv\"os University, Budapest, Hungary.
\newline
{\sf E-mail: bezdek@math.ucalgary.ca}

\end{document}